\newcount\strangepages \strangepages = 0
\def\lastchange{November 30, 2006}
\ifnum \strangepages = 1 \voffset -100 pt  \fi
\newcount\listcodes
\ifnum \listcodes > 0 {
\else  \fi
\ifnum \listcodes > 0 {\def\labeleqn#1{
\label{#1} {\quad\kern-6em\mbox{\tt \{#1\}}\quad}  }
\else  \def\labeleqn#1{\label{#1}}  \fi

\def\hdq{\hskip 0.25 em\relax}

\documentclass[12pt]{article}

\ifnum \strangepages = 1 \textheight = 6 in \fi

\renewcommand{\labelenumi}{{\normalfont (\alph{enumi})}}

\usepackage{amsmath}
\usepackage{amssymb}
\usepackage{enumerate}

\def\x#1{\def\@@eqncr{\let\@tempa\relax
  \ifcase\@eqcnt \def\@tempa{& & &}\or \def\@tempa{& &}
     \else \def\@tempa{&}\fi & #1\cr}}

\font\teneufm=eufm10
\font\seveneufm=eufm7
\font\fiveeufm=eufm5
\newfam\eufmfam
\textfont\eufmfam=\teneufm
\scriptfont\eufmfam=\seveneufm
\scriptscriptfont\eufmfam=\fiveeufm

\newfam\msbfam
\font\tenmsb=msbm10 scaled \magstep1 \textfont\msbfam=\tenmsb
\font\sevenmsb=msbm7 scaled \magstep1  \scriptfont\msbfam=\sevenmsb
\font\fivemsb=msbm5 scaled \magstep1  \scriptscriptfont\msbfam=\fivemsb
\def\Bbb{\fam\msbfam \tenmsb}

\newfam\msbfamx
\font\tenmsbx=msbm10  \textfont\msbfamx=\tenmsbx
\font\sevenmsbx=msbm7   \scriptfont\msbfamx=\sevenmsbx
\font\fivemsbx=msbm5  \scriptscriptfont\msbfamx=\fivemsbx

\def\RR{{\Bbb R}}

\def\ra{\rightarrow}
\def\sgn{\hbox{\rm sgn}\,}
\def\ss{\subseteq}

\def\tombstone{\rule{.08in}{.1in}}
\def\endpf{\hfill\tombstone\medskip}
\def\endpfarray{\vskip - 4 ex \rightline{\tombstone}\medskip}

\def\littleo#1{o(x^{#1})}
\def\bigo#1{O(x^{#1})}

\newtheorem{theorem}{Theorem}[section]
\newtheorem{observation}[theorem]{Observation}
\newtheorem{lemma}[theorem]{Lemma}
\newtheorem{proposition}[theorem]{Proposition}
\newtheorem{corollary}[theorem]{Corollary}

\newtheorem{example}[theorem]{Example}
\newtheorem{remark}[theorem]{Remark}
\newtheorem{definition}[theorem]{Definition}
\newtheorem{notation}[theorem]{Notation}
\def\btheorem{\begin{theorem}\sl}
\def\btheoremc#1{\begin{theorem}[#1]\sl}
\def\etheorem{\end{theorem}}
\def\blemma{\begin{lemma}\sl}
\def\blemmac#1{\begin{lemma}[#1]\sl}
\def\elemma{\end{lemma}}
\def\bproposition{\begin{proposition}\sl}
\def\bpropositionc#1{\begin{proposition}[#1]\sl}
\def\eproposition{\end{proposition}}
\def\bcorollary{\begin{corollary}\sl}
\def\bcorollaryc#1{\begin{corollary}[#1]\sl}
\def\ecorollary{\end{corollary}}
\def\bexample{\begin{example}\rm}
\def\bexamplec#1{\begin{example}[#1]\rm}
\def\eexample{\endpf\end{example}}
\def\bremark{\begin{remark}\rm}
\def\bremarkc#1{\begin{remark}[#1]\rm}
\def\eremark{\end{remark}}
\def\bdefinition{\begin{definition}\rm}
\def\bdefinitionc#1{\begin{definition}[#1]\rm}
\def\edefinition{\end{definition}}
\def\bnotation{\begin{notation}\rm}
\def\bnotationc#1{\begin{notation}[#1]\rm}
\def\enotation{\end{notation}}
\def\bobservation{\begin{observation}\rm}
\def\bobservationc#1{\begin{observation}[#1]\rm}
\def\eobservation{\end{observation}}
\def\bequation{\begin{equation}}
\def\eequation{\end{equation}}
\def\bprop{\begin{proposition}\sl}
\def\bpropc#1{\begin{proposition}[#1]\sl}

\begin{document}
\quad

\vskip - 1.375in
\rightline{\footnotesize\lastchange}

\vskip 1 in
\begin{center}
{\LARGE \bf The Lagrange Inversion Theorem}
\smallskip \\
{\LARGE \bf in the Smooth Case}\footnote{The authors are happy to
thank the American Institute of Mathematics for its hospitality during
the writing of this paper.}
\bigskip \\
by Harold R. Parks and Steven G. Krantz
\end{center}

\begin{quote} {\bf Abstract:} \sl
The classical Lagrange inversion theorem is a concrete,
explicit form of the implicit function theorem for
real analytic functions.  The authors derive a suitable
version of this result for $C^\infty$ functions.  Along
the way, they find a new asymptotic for smooth functions.
\end{quote}

\def\fin{\end{document}}


\section{The Problem}

The implicit function theorem has a long and colorful history.  Finding
its provenance in considerations of problems of celestial mechanics
(as studied by Lagrange and Cauchy, among others), the result was
at first a rather primitive observation about monotone functions
on $\RR^1$.  Over time, the result was extended to $N$ variables, and
the monotonicity hypothesis was replaced by the now more familiar
assumption of nondegeneracy of the Jacobian at a point (see [KPb] for
a more detailed history). 

It is easy to imagine that there were a number of vestigial forms
of the implicit function theorem that historically preceded
the crisp result that can be found in textbooks today.	One of
these is the so-called 
{\it Lagrange inversion theorem}.\footnote{Lagrange's inversion formula 
can be found in [LAG]. 
%
%
The  date of the memoir is 1768, but it was
published in  1770 in
{\it Histoire de L'Acad\'emie Royal des Sciences et Belles-Lettres}.}  That
is the topic of the present paper.  
The classical Lagrange theorem
is about 
%
%
%
analytic functions.
Our purpose here is to extend the result to
$C^k$ or $C^\infty$ functions.

The form of the Lagrange inversion theorem that we will consider
in this paper concerns solving
\begin{equation}\label{eq:main}
y = x + f(y)
\end{equation}
for $y$ as a function of $x$. Lagrange's result is that, 
if $f(0)=0$, $|f'(0)|<1$, and $f$ is real analytic, then
\begin{eqnarray}
\nonumber
y &=& x + \sum_{n=1}^\infty \frac1{n!}
\ \left(\frac{d\ }{dx}\right)^{n-1}
\left\{f^n(x)\right\}
\\
\label{eq:lagrange.formula}
&=&
x + f(x) + \sum_{n=1}^\infty \frac1{(n+1)!}
\ \left(\frac{d\ }{dx}\right)^{n}
\left\{f^{n+1}(x)\right\}
\,.
\end{eqnarray}
One can obtain (\ref{eq:lagrange.formula}) from  
[KPb; Theorem~2.3.1] by complexifying $f$,
substituting $\phi(z) = f(z)$, and taking $\psi(z) \equiv z$, $t=1$.
A 
version of the Lagrange theorem in the context of real analysis, for smooth
functions, is to be found in [GRO].  The arguments presented there seem to
be incomplete, and the statement of the result incorrect.  The purpose of
the present paper is to explore these matters further and to set the
record straight.

We will find it useful to combine the Lagrange inversion formula
with the formula of Fa\'a di Bruno.\footnote{Fa\'a di Bruno's
formula first appeared in [FDB]---for a proof, see Section~1.3 of
[KPa].}
The formula of Fa\'a di Bruno tells us that if $I$ and $J$ are 
open intervals in $\RR$ and
$f:I\to J$ and $h:J\to\RR$ are $C^\infty$ functions,
then the $n$th
derivative of $ h \circ f$  is given by
\begin{eqnarray}
&&
\kern - 3 em 
(h\circ f)^{(n)}(t) =  \nonumber\\
&&
\kern - 2 em 
\sum\frac{n!}{k_1!\,k_2!\dots k_n!}\ h^{(k)} \bigl (f(t) \bigr )
\left(\frac{f^{(1)}(t)}{1!}\right)^{k_1}
\left(\frac{f^{(2)}(t)}{2!}\right)^{k_2}
\dots\left(\frac{f^{(n)}(t)}{n!}\right)^{k_n},
\label{eq:faa.de.bruno}
\end{eqnarray}
where $k=k_1+k_2+\cdots+k_n$ and the sum is taken over all
$k_1,k_2,\dots,k_n$
for which $k_1 + 2k_2 +\cdots +nk_n = n.$

In particular, when $h(t) = t^{n+1}$,
we have 
$$
h^{(k)}(t) = 
\left\{\begin{array}{cl}
[\,(n+1)!/(n+1-k)!\,]\,t^{n+1-k} &\mbox{if\ }k\leq n+1,\\[1ex]
0 &\mbox{if\ }n+1<k.
\end{array}\right.
$$
Fa\'a di Bruno's formula then becomes
\begin{eqnarray}
&&
\kern - 3 em 
\left(\frac{d\ }{dx}\right)^{n}
\left\{f^{n+1}(x)\right\}  (t) =  \nonumber\\
&&
\kern - 2 em 
\sum\frac{n! \,(n+1)!}{k_0!\,k_1!
\dots k_n!}
\left(\frac{f^{(0)}(t)}{0!}\right)^{k_0}
\left(\frac{f^{(1)}(t)}{1!}\right)^{k_1}
\dots\left(\frac{f^{(n)}(t)}{n!}\right)^{k_n} ,
\label{eq:faa.de.bruno2}
\end{eqnarray}
where $k_0=(n+1)-(k_1+k_2+\cdots+k_n)$ and the sum is taken over all
$k_1,k_2,\dots,k_n$
for which $k_1 + 2k_2 +\cdots +nk_n = n.$
Accordingly, Lagrange's inversion formula can be written
\begin{equation}\label{eq:lagrange.formula.bruno}
y = x + f(x) + \sum_{n=1}^\infty n!
\sum
\prod_{i=0}^n \frac{1}{k_i!}
\left(\frac{f^{(i)}(x)}{i!}\right)^{k_i}
\,,
\end{equation}
where the inner sum is taken over all
$k_1,k_2,\dots,k_n$
for which $k_1 + 2k_2 +\cdots +nk_n = n$ and where 
$k_0$ is defined by $k_0=(n+1)-(k_1+k_2+\cdots+k_n)$.
%
%

We can rewrite (\ref{eq:lagrange.formula.bruno}) as
\begin{equation}\label{eq:lagrange.formula.bruno1}
y = x + f(x) + 
\sum_{(k_0,\dots,k_m)\in \Xi}
C_{k_0,\dots,k_m}\cdot
\left[\,f(x)\,\right]^{k_0}
\cdot
\prod_{i=1}^m 
\left[\,f^{(i)}(x)\,\right]^{k_i}
\,,
\end{equation}
where 
%
%
\begin{eqnarray*}
\Xi = \bigcup_{n=1}^\infty
\bigcup_{m=1}^n
\Big\{\,(k_0,k_1,\dots,k_m) &:& 0\leq k_i,\mbox{\ for\ } i=1,\dots,m-1,\\
&&
\ 1\leq k_m,\\
&&
\ k_0 = (n+1)-(k_1+k_2+\cdots+k_m),\\
&&
\ k_1 + 2k_2 +\cdots +mk_m = n\,\Big\}
\,,
\end{eqnarray*}
and where
$$
C_{k_0,\dots,k_m} =
\frac{(k_1+2k_2+\cdots+mk_m)! }{k_0!}
\,\left(\prod_{i=1}^m k_i!\,(i!)^{k_i}\right)^{-1}
\,,
$$
for $(k_0,\dots,k_m)\in\Xi$. Note that all the 
$C_{k_0,\dots,k_m}$ are positive.

%

\section{\boldmath A Counterexample in the $C^\infty$ Case}
We will now show that the Lagrange inversion formula
is not true in the $C^\infty$ category.

To begin with, we let $f$ be a real analytic function 
defined on an open interval containing $[0,1]$  and
such that
\def\labelenumi{{\rm(\alph{enumi})}}
\def\theenumi{\alph{enumi}}
\begin{enumerate}
\item
$f(0) = 0$, 
\item\label{item:example.condition}
$0< f'(0)<1$,
\item
for each non-negative integer $m$,
$0< f^{(m)} (t) $ holds for all $0<t$.
\end{enumerate}
Note that, by (\ref{item:example.condition}), we have
\begin{equation}\label{eq:example.condition}
0<f(t)<t \hbox{\ \ for \ \ }t\in (0,1)
\,.
\end{equation}

Many such functions exist. For example, we could 
take $f(x)= \lambda\,e^x - \lambda$, where
$0<\lambda<1$. 

Since $f$ is real analytic, the Lagrange
inversion formula is valid in an open interval $I$ containing 0.

\medskip
Next, we will define sequences $\{a_\ell\}$ and $\{b_\ell\}$ in $I$
that  satisfy
\begin{eqnarray*}
&& b_0>a_0>b_1>a_1>\cdots b_\ell>a_\ell > \cdots 0
\\
&&
b_\ell = a_\ell + f(b_\ell)
\,.
\end{eqnarray*}
In fact, we may choose any  $b_0$ in $I$ with 
$0<b_0<1$.  We set
$a_0 = b_0 - f(b_0)$, noting that (\ref{eq:example.condition}) implies that
$0<a_0<b_0$. Proceeding inductively, we observe that
if $ b_0>a_0>b_1>a_1>\cdots b_\ell>a_\ell $ have already 
been chosen, then we may choose any $b_{\ell+1}$ with
$0< b_{\ell+1} < a_\ell$ and set 
$$
a_{\ell+1} = b_{\ell+1}  - f(b_{\ell+1})
\,.
$$
Again we use (\ref{eq:example.condition}) to conclude that
$0<a_{\ell+1}< b_{\ell+1} $.

The Lagrange inversion formula (\ref{eq:lagrange.formula.bruno1}) tells us that
\begin{equation}\label{eq:example1}
b_\ell = a_\ell + f(a_\ell) + 
\sum_{(k_0,\dots,k_m)\in \Xi}
C_{k_0,\dots,k_m}\cdot
\left[\,f(a_\ell)\,\right]^{k_0}
\cdot
\prod_{i=1}^m 
\left[\,f^{(i)}(a_\ell)\,\right]^{k_i}
\end{equation}
holds, for $\ell = 0,1,\dots$.

\medskip
Now we consider a $C^\infty$ function $g:\RR\to\RR$
such that
\begin{enumerate}\setcounter{enumi}{3}
\item \label{item:property_g1}
$g(b_\ell) = 0$ for $\ell = 0,1,\dots$,
\item \label{item:property_g2}
$g(a_\ell) > 0$ for $\ell = 0,1,\dots$,
\item \label{item:property_g3}
$g^{(m)}(a_\ell) \geq 0$ for $\ell = 0,1,\dots$
and $m = 1,2, \dots$.
\end{enumerate}
[We consider in Lemma~\ref{lemma:apply.borel} 
and Remark~\ref{remark:apply.lemma} why such a function $g$ exists.]
Then we have
\begin{equation}\label{eq:c_infinity_places}
b_\ell = a_\ell + (f+g)(b_\ell)
\,.
\end{equation}
If the Lagrange inversion formula held for the function 
$f+g$ in any neighborhood
of $0$, then it would hold
at infinitely many of the pairs $(a_\ell, b_\ell)$ in
(\ref{eq:c_infinity_places}).   
For such a  pair $(a_\ell, b_\ell)$,
 the Lagrange inversion formula 
for the function $f+g$ would tell us that 
\begin{equation}\label{eq:example2}
b_\ell = a_\ell + (f+g)(a_\ell) + 
\sum_{(k_0,\dots,k_m)\in \Xi}
C_{k_0,\dots,k_m}\cdot
\left[\,(f+g)(a_\ell)\,\right]^{k_0}
\cdot
\prod_{i=1}^m 
\left[\,(f+g)^{(i)}(a_\ell)\,\right]^{k_i}
\end{equation}
holds. But, in fact,
 (\ref{eq:example2}) fails to hold for any choice of $\ell$.  This is so because, while the 
lefthand sides of (\ref{eq:example1}) and (\ref{eq:example2})
are equal, 
 the righthand side of 
(\ref{eq:example2}) consists of  the righthand side of
(\ref{eq:example1}) plus nonnegative  terms, among which 
is the positive term  $g(a_\ell)$.

\medskip
Finally,
it remains to show that a function $g$ satisfying the above conditions
 exists, but this fact follows from the next lemma.
\begin{lemma}  \label{lemma:apply.borel}
Let $\{x_n\}_{n=1}^\infty$ be a strictly
decreasing sequence of positive numbers with $\lim_{n\ra\infty}x_n=0$.
If, for each $n=1,2,\dots,$ there is given a sequence
$\{\sigma_{n,k}\}_{k=0}^\infty \ss \{\,-1,\,0,\,1\,\}$, then
there exists a $C^\infty$ function $g:\RR\to\RR$ such that
\begin{eqnarray}
\label{eq:property_gk}
\sgn\Big[\, g^{(k)}(x_n)\,\Big] &=& \sigma_{n,k}\,,\mbox{\ for $n=1,2,\dots$ and $k=0,1,\dots$,} \\
\label{eq:property_g_at_0}
g^{(k)}(0) &=& 0\,,\mbox{\ for $n=1,2,\dots$.}
\end{eqnarray}
\end{lemma}
{\bf Proof.}  Let $I_1, I_2,\dots$ be pairwise disjoint open intervals with
$x_n \in I_n$. Assume also that $I_1$ is bounded. The other intervals are
automatically bounded, because $I_n\subset (x_{n+1},\,x_{n-1})$
holds for $n=2,3,\dots$.

An old theorem of \'Emile Borel (see page 44 of [BOR] or else [KPa])
tells us that, for each $n$, there exists a $C^\infty$ function $\phi_n$ with
\begin{equation}\label{eq:choice_of_phi}
\phi^{(k)}_n (x_n) = \sigma_{n,k}\,,\mbox{\ for $k=0,1,\dots$.}
\end{equation}
We may assume that $\phi_n$ has compact support contained in $I_n$.

We will define
\begin{equation}\label{eq:def_of_g}
g(x) = \sum_{k=1}^\infty \xi_n^{-1}\, \phi_n(x)
\,,
\end{equation}
where the positive numbers $\xi_n$ will be chosen large enough
to insure that $g$ is $C^\infty$.  Specifically, it is
easy to see that the choices
$$
\xi_n = e^{1/x_{n+1}}\cdot 
\sup\left\{\,\left|\phi_m^{(k)}(x)\right| : 
x\in I_m,\ 1\leq m\leq n,\ 0\leq k\leq n\,\right\}
$$ 
will suffice to guarantee that $g$ is $C^\infty$
and that (\ref{eq:property_g_at_0}) holds.

It is clear that (\ref{eq:choice_of_phi}) 
and  (\ref{eq:def_of_g}) insure that (\ref{eq:property_gk})
holds.
\endpf

\bremark \label{remark:apply.lemma}
Lemma~\ref{lemma:apply.borel} is applied to construct 
a function $g$ satisfying (\ref{item:property_g1}),
(\ref{item:property_g2}), and (\ref{item:property_g3})
by setting 
$$
x_{2\ell+1} = b_\ell\,,\quad
\sigma_{2\ell+1,0} =0\,,\quad
x_{2\ell+2} = a_\ell\,,\quad
\sigma_{2\ell+2,0} = 1\,,
\hbox{\ \ for\ \ }\ell=0,1,2,\dots
$$
and setting
$$
\sigma_{n,k} = 0\hbox{\ \ for\ \ }n=1,2,\dots
\hbox{\ \ and\ \ }k=1,2,\dots\,.
$$
\eremark

\section{\boldmath A Positive Result in the $C^N$ Case}

It is reasonable to conjecture---and we will 
show it to be true---that, if $f$ is 
$C^N$, then 
the function $y$ satisfying $y=x+f(y)$ is
well approximated by the {\em truncated Lagrange formula of order $N$}.
By the later we mean that part of 
$$
x + \sum_{n=1}^\infty \frac1{n!}
\ \left(\frac{d\ }{dx}\right)^{n-1}
\left\{f^n(x)\right\}
$$
that only contains derivatives of $f$ of order $N$ or 
smaller (the truncated Lagrange formula is given below in
(\ref{eq:lagrange.formula.brunoN})).

We need a lemma.

\blemma\label{lemma:taylor}
Suppose $f(x)$ is $N$ times continuously differentiable,
$N\geq 1$, in a 
neighborhood of $x=0$ with $f(0)= 0$ and $|f'(0)|< 1$.
Let $P(x)$ be the Taylor polynomial of degree $N$ for
$f$ at $0$.
If $y=y(x)$ satisfies 
$$
y = x + f(y)
$$
in a neighborhood of $0$, then
$$
y - \left[
x + \sum_{n=1}^\infty \frac1{n!}
\ \left(\frac{d\ }{dx}\right)^{n-1}
\left\{P^n(x)\right\} 
\right]
= \littleo{N}
$$
as $x\rightarrow 0$. 
\elemma
\bremark 
Because $P(x)$ is the Taylor polynomial of degree $N$ for
$f$ at $0$, only derivatives of 
$f$ of order $N$ or less are involved in the construction of 
$P$. Thus we see that
the expression
$$
 x + \sum_{n=1}^\infty \frac1{n!}
\ \left(\frac{d\ }{dx}\right)^{n-1}
\left\{P^n(x)\right\}
$$
appearing in the lemma involves only derivatives of 
$f$ of order $N$ or less, 
but note that those derivatives are 
evaluated at $0$ rather than at $x$ as are the derivatives   
in Lagrange's formula.
\eremark
{\bf Proof of the Lemma.\  }
By Taylor's theorem, we may write
\begin{equation}\label{eq:taylor}
f(x) = P(x) + R(x)
\,,
\end{equation}
where
\begin{equation}\label{eq:behavior.R}
R(x) = \littleo{N} \hbox{\  as\ }x\rightarrow 0
\,.
\end{equation}

Since
$f$ is at least $C^1$ near $0$, we can apply the 
implicit function theorem to see that there 
exists a function
$y = y(x)$ satisfying
\begin{equation}\label{eq:basic.lagrange.update}
y = x + f(y)
\,.
\end{equation}
Implicit differentiation shows us that
$y'(0) = 1/[1-f'(0)]$. Thus   
we have  
\begin{equation}\label{eq:behavior.y}
y = \bigo{} \hbox{\  as\ } x\rightarrow 0
\,.
\end{equation} 
Combining (\ref{eq:taylor}) and (\ref{eq:basic.lagrange.update}),
we have
\begin{equation}\label{eq:taylor.lagrange}
y = x + P(y) + R(y)
\,.
\end{equation}

Since $P(x)$ is an analytic function,
we  may apply the Lagrange inversion theorem to see that 
the  function $z = z(x)$ 
given by
\begin{equation}\label{eq:lagrange.formula.z}
z = x + \sum_{n=1}^\infty \frac1{n!}
\ \left(\frac{d\ }{dx}\right)^{n-1}
\left\{P^n(x)\right\}
\end{equation}
satisfies
\begin{equation}\label{eq:lagrange.z}
z = x + P(z)
\,.
\end{equation}
As before, we have 
\begin{equation}\label{eq:behavior.z}
z = \bigo{} \hbox{\  as\ } x\rightarrow 0
\,.
\end{equation} 

Note that (\ref{eq:lagrange.formula.z}) gives us
\begin{equation}\label{eq:target}
y-\left(\,
x + \sum_{n=1}^\infty \frac1{n!}
\ \left(\frac{d\ }{dx}\right)^{n-1}
\left\{P^n(x)\right\}
\,\right) = y-z
\,.
\end{equation}

Now consider an $x$ near $0$ for which $y(x)\neq z(x)$.
Applying the mean value theorem to 
$$
Q(t) = t - P(t)
$$
on the interval $[\min\{y,z\},\hdq \max\{y,z\}]$,
we obtain
 $\xi$ 
in that interval
with
\begin{eqnarray}
\nonumber
(y-z) - [P(y)-P(z)]
&=&
Q(y)-Q(z)\\[1.5ex]
\nonumber
&=& Q'(\xi) (y-z)\\[1.5ex]
\label{eq:mean.value}
&=& [1-P'(\xi)]\,(y-z)
\,.
\end{eqnarray}
By (\ref{eq:taylor.lagrange}),
(\ref{eq:lagrange.z}),
and (\ref{eq:mean.value}),  we have
\begin{eqnarray}
\nonumber
R(y) &=& [y - x - P(y)]
-[z-x-P(z)] \\[1.5ex]
\nonumber
&=& (y-z) -[P(y)-P(z)]\\[1.5ex]
\label{eq:combine.mean}
&=& [1-P'(\xi)]\,(y-z)
\,.
\end{eqnarray}
By (\ref{eq:behavior.R}) and (\ref{eq:behavior.y}), we 
have $R(y) = \littleo{N}$, so we see from (\ref{eq:combine.mean}) that 
\begin{equation}\label{eq:penultimate.step}
[1-P'(\xi)]\,(y-z) = \littleo{N}
\,.
\end{equation}
Since $P'(0) = f'(0) < 1$ and since $\xi\rightarrow 0$ as
$x\rightarrow 0$ [by (\ref{eq:behavior.y})
and (\ref{eq:behavior.z})], we conclude from 
(\ref{eq:penultimate.step}) that
$y-z = \littleo{N}$. The result follows from (\ref{eq:target}).
\endpf

\btheorem
Suppose $f(x)$ is $N$ times continuously differentiable,
$N\geq 1$, in a 
neighborhood of $x=0$ with $f(0)= 0$ and $|f'(0)|< 1$.
If $y=y(x)$ satisfies 
$$
y = x + f(y)
$$
in a neighborhood of $0$, then
$$
y - \left[
x + f(x) + 
\sum_{(k_0,\dots,k_m)\in \Xi_N}
C_{k_0,\dots,k_m}\cdot
\left[\,f(x)\,\right]^{k_0}
\cdot
\prod_{i=1}^m 
\left[\,f^{(i)}(x)\,\right]^{k_i}
\right]
= \littleo{N}
$$
as $x\rightarrow 0$. 
Here 
%
%
\begin{eqnarray*}
\Xi_N = \bigcup_{n=1}^\infty
\bigcup_{m=1}^{\min\{n,N\}}
\Big\{\,(k_0,k_1,\dots,k_m) &:& 0\leq k_i,\mbox{\ for\ } i=1,\dots,m-1,\\
&&
\ 1\leq k_m,\\
&&
\ k_0 = (n+1)-(k_1+k_2+\cdots+k_m),\\
&&
\ k_1 + 2k_2 +\cdots +mk_m = n\,\Big\}
\,,
\end{eqnarray*}
and 
%
%
$$
C_{k_0,\dots,k_m} =
\frac{(k_1+2k_2+\cdots+mk_m)! }{k_0!}
\,\left(\prod_{i=1}^m k_i!\,(i!)^{k_i}\right)^{-1}
\,.
$$
\etheorem
\bremark
By (\ref{eq:lagrange.formula.bruno1}), we see that
the expression
\begin{equation}\label{eq:lagrange.formula.brunoN}
x + f(x) + 
\sum_{(k_0,\dots,k_m)\in \Xi_N}
C_{k_0,\dots,k_m}\cdot
\left[\,f(x)\,\right]^{k_0}
\cdot
\prod_{i=1}^m 
\left[\,f^{(i)}(x)\,\right]^{k_i}
\,,
\end{equation}
appearing in the theorem  involves only derivatives of 
$f$ of order $N$ or smaller, and those derivatives are 
evaluated
 at $x$ as are the derivatives   in Lagrange's formula. Thus 
 (\ref{eq:lagrange.formula.brunoN})
is the truncated Lagrange formula  described at the beginning of this 
section.
\eremark
{\bf Proof of the Theorem.\ }
As in Lemma~\ref{lemma:taylor},
let $P(x)$ be the Taylor polynomial of degree $N$ for $f$
at $0$.
The formula of Fa\'a di Bruno tells us that
\begin{eqnarray*}
&&
x + \sum_{n=1}^\infty \frac1{n!}
\ \left(\frac{d\ }{dx}\right)^{n-1}
\left\{P^n(x)\right\} 
\\[1ex]
&&
=
\ x + P(x) + 
\sum_{(k_0,\dots,k_m)\in \Xi}
C_{k_0,\dots,k_m}\cdot
\left[\,P(x)\,\right]^{k_0}
\cdot
\prod_{i=1}^m 
\left[\, P^{(i)}(x) \,\right]^{k_i}
\,.
\end{eqnarray*}
Since $P$ is a polynomial of degree $N$,
we have $P^{(i)}(x)\equiv 0$ whenever $i>N$. Thus
we have
\begin{eqnarray*}
&&
\ x + P(x) + 
\sum_{(k_0,\dots,k_m)\in \Xi}
C_{k_0,\dots,k_m}\cdot
\left[\,P(x)\,\right]^{k_0}
\cdot
\prod_{i=1}^m 
\left[\, P^{(i)}(x) \,\right]^{k_i}
\\[1ex]
&&
=
\ x + P(x) + 
\sum_{(k_0,\dots,k_m)\in \Xi_N}
C_{k_0,\dots,k_m}\cdot
\left[\,P(x)\,\right]^{k_0}
\cdot
\prod_{i=1}^m 
\left[\, P^{(i)}(x) \,\right]^{k_i}
\,.
\end{eqnarray*}

Applying Lemma~\ref{lemma:taylor},
we see that
\begin{eqnarray*}
&& 
\littleo{N} = y-
\left[
x + \sum_{n=1}^\infty \frac1{n!}
\ \left(\frac{d\ }{dx}\right)^{n-1}
\left\{P^n(x)\right\} 
\right]
\\[1ex]
&&\quad
=
\ y - \left[
x + P(x) + 
\sum_{(k_0,\dots,k_m)\in \Xi_N}
C_{k_0,\dots,k_m}\cdot
\left[\,P(x)\,\right]^{k_0}
\cdot
\prod_{i=1}^m 
\left[\,P^{(i)}(x)\,\right]^{k_i}
\right]
\\[1ex]
&&\quad
=
\ y - \left[
x + f(x) + 
\sum_{(k_0,\dots,k_m)\in \Xi_N}
C_{k_0,\dots,k_m}\cdot
\left[\,f(x)\,\right]^{k_0}
\cdot
\prod_{i=1}^m 
\left[\,f^{(i)}(x)\,\right]^{k_i}
\right]
\\[1ex]
&&\qquad+\ [f(x) -P(x)]
\\[1ex]
&&\qquad+ \ \sum_{ (k_0,\dots,k_m)\in \Xi_N }
C_{k_0,\dots,k_m}
\\[1ex]
&& \qquad\qquad \cdot
\ \left(
\left[\,f(x)\,\right]^{k_0}
\cdot
\prod_{i=1}^m 
\left[\,f^{(i)}(x)\,\right]^{k_i}
-
 \left[\, P(x)\,\right]^{k_0}
\cdot
\prod_{i=1}^m 
\left[\, P^{(i)}(x) \,\right]^{k_i}
\right)
\end{eqnarray*}

Of course, we have $f(x)-P(x)= \littleo{N} $. So it will suffice
to show that
\begin{equation}\label{eq:key.to.bound}
\left[\,f(x)\,\right]^{k_0}\cdot
\prod_{i=1}^m 
\left[\,f^{(i)}(x)\,\right]^{k_i} 
-
\left[\,P(x)\,\right]^{k_0}\cdot
\prod_{i=1}^m 
\left[\, P^{(i)}(x) \,\right]^{k_i}
= \littleo{N}
\end{equation}
holds, for each $ (k_0,\dots,k_m)\in \Xi_N $.

We rewrite the lefthand side of (\ref{eq:key.to.bound}) as 
\begin{eqnarray*}
&& \left[\,f(x)\,\right]^{k_0}\cdot
\prod_{i=1}^m
\left[\,f^{(i)}(x)\,\right]^{k_i} 
-
\left[\,f(x)\,\right]^{k_0} \cdot
\prod_{i=1}^{m-1} 
\left[\,f^{(i)}(x)\,\right]^{k_i} 
\cdot
\left[\, P^{(m)}(x) \,\right]^{k_m}\\[1ex]
&&+\cdots+\\[1ex]
&&\qquad+ \bigg(\left[\,f(x)\,\right]^{k_0}\cdot
\prod_{i=1}^{j} 
\left[\,f^{(i)}(x)\,\right]^{k_i} 
\cdot \prod_{i=j+1}^m 
\left[\, P^{(i)}(x) \,\right]^{k_i}\\[1ex]
&&\qquad\qquad\qquad-
 \left[\,f(x)\,\right]^{k_0}\cdot
\prod_{i=1}^{j-1} 
\left[\,f^{(i)}(x)\,\right]^{k_i} 
\cdot
\prod_{i=j}^m 
\left[\, P^{(i)}(x) \,\right]^{k_i}\bigg) \\[1ex]
&&+\cdots+\\[1ex]
&&
+\ \left[\,f(x)\,\right]^{k_0}
\prod_{i=1}^m
\left[\, P^{(i)}(x) \,\right]^{k_i}
-
\left[\,P(x)\,\right]^{k_0}
\prod_{i=1}^m 
\left[\, P^{(i)}(x) \,\right]^{k_i}
\,.
\end{eqnarray*}
Thus it will suffice to show that
$$
\left[\,f(x)\,\right]^{k_0}\,\bigg(
\left[\, P^{(j)}(x)\,\right ]^ {k_j} - \left[\,f^{(j)}(x)\,\right]^{k_j}
\bigg)
=\littleo{N}
$$
holds, for each $j=1,2,\dots,m$ for which $k_j>0$, and that
$$
\left[\,f(x)\,\right]^{k_0}
-
\left[\,P(x)\,\right]^{k_0}
=\littleo{N}
\,.
$$

If $j\in\{1,2,\dots,m\}$ and $k_j>0$ holds, then we have
\begin{eqnarray*}
&&\left[\, P^{(j)}(x)\,\right ]^ {k_j} - \left[\,f^{(j)}(x)\,\right]^{k_j}
\\[1ex]
&&\qquad =
\left[\, P^{(j)}(x) - f^{(j)}(x) \,\right] \cdot \sum_{\ell=0}^{k_j-1}
\left[\, P^{(j)}(x) \,\right]^ {\ell} \left[ \,f^{(j)}(x) \,\right]^ {k_j-1-\ell}
=
\littleo{N-j}
\,.
\end{eqnarray*}
Using $f(x)= \bigo{1}$ and $P^{(j)}(x) - f^{(j)}(x) = \littleo{N-j}$, 
we obtain
$$
\left[\,f(x)\,\right]^{k_0}\,\bigg(
\left[\, P^{(j)}(x)\,\right ]^ {k_j} - \left[\,f^{(j)}(x)\,\right]^{k_j}
\bigg)
=
\littleo{N-j+ k_0}
\,.
$$
By the definition of $k_0$ and using
$ n = k_1+ 2 k_2+\cdots+k_j+\cdots+mk_m $, we have
\begin{eqnarray*}
&& N-j+k_0
\\[1ex]
&&\qquad 
= \ N-j + (n+1)-(k_1+  k_2+\cdots+k_j+\cdots+k_m)\\[1ex]
&&\qquad
=\ N -j + 1 + k_2 + 2k_3+ \cdots (j-1)k_j + \cdots +(m-1)k_m\\[1ex]
&&\qquad 
=\ N  + (j-1)(k_j-1) + 
\sum_{ \stackrel{\ell=2,\dots,m}
 {\ell\neq j}}(\ell-1)k_\ell
\geq N
\,.
\end{eqnarray*}

Finally, we observe that
\begin{eqnarray*}
\left[\,f(x)\,\right]^{k_0}
-
\left[\,P(x)\,\right]^{k_0}
&=&
\Big(\, f(x) - P(x) \,\Big)
\sum_{i=1}^{k_0-1} 
\left[\,f(x)\,\right]^{i}
\left[\,P(x)\,\right]^{k_0-i-1}\\[1ex]
&=& \littleo{N}
\,.
\end{eqnarray*}
\endpfarray

\noindent {\Large \sc References}

\begin{enumerate}
\item[]
\begin{enumerate}
\itemindent= 0 em
\labelsep = 0.5 em
\labelwidth = 0 em
\item[{\bf [BOR]}]  \'Emile Borel, Sur quelques points de la
th\'eorie des fonctions, {\it Annales Scientifiques de l'\'Ecole Normale
Sup\'erieure} (3) 12 (1895), 9--55
\item[{\bf [FDB]}]  Francesco Fa\'a di Bruno, Note sur une nouvelle formule de calcul 
diff{\'e}rentiel. {\it Quarterly Journal of Pure and Applied Mathematics} 
1 (1857), 359--360.
\item[{\bf [GRO]}]  Nathaniel Grossman, A $C^\infty$ Lagrange inversion theorem,
{\it American  Mathematical  Monthly} 112 (2005), 512--514.
\item[{\bf [KPa]}]  Steven G. Krantz \& Harold R. Parks, {\it A Primer of
Real Analytic Functions}, second edition,  Birkh\"{a}user, Boston, 2002.
\item[{\bf [KPb]}]  \vrule height 0.5 pt depth 0 pt width 3 em 
\kern 0.5 em  {\it The Implicit Function
Theorem}, Birkh\"{a}user, Boston, 2002.
\item[{\bf [LAG]}] Joseph Louis Lagrange,
Nouvelle m\'ethode pour r\'esoudre les \'equations litt\'erales par le 
moyen des s\'eries,
({\it M\'emoires de l'Acad\'emie Royale
des Sciences et Belles-Lettres de 
Berlin} {\bf 24})
{\it {\OE}uvres de Lagrange}, volume 3,
Gauthier-Villars, Paris, 1869, p. 25.
\end{enumerate}
\end{enumerate}

\fin